\documentclass[11pt,leqno]{amsart}
\usepackage{amsmath}
\usepackage{amsthm}
\usepackage{amssymb}
\theoremstyle{definition}

\theoremstyle{plain}

\begin{document}
\setcounter{page}{63}
\begin{flushleft}
\scriptsize{{Journal of the Nigerian Association of Mathematical Physics\/} {\bf 11}\,(November 2007), 
\\
Pages 63--66}
\end{flushleft}
\vspace{16mm}

\begin{center}
{\normalsize\bf On some functionals associated with certain coefficient problems} \\[12mm]
     {\sc K. O. BABALOLA$^{1}$} \\ [8mm]

\begin{minipage}{123mm}
{\small {\sc Abstract.}
Under certain conditions, we obtain sharp bounds on some functionals defined in the coefficient space of starlike functions. It has been found that the functionals are closely associated with certain coefficient problems, which are of independent interest.}
\end{minipage}
\end{center}

 \renewcommand{\thefootnote}{}
 \footnotetext{$^1$\,Department of Mathematics, University of Ilorin, Ilorin, Nigeria. E-mail:
ummusalamah.kob\symbol{64}unilorin.edu.ng}

\def\iff{if and only if }
\def\S{Smarandache }
\newcommand{\norm}[1]{\left\Vert#1\right\Vert}

\vskip 12mm

{\bf 1. Introduction }
\medskip

Let $S$ denote the family of functions:
$$f(z)=z+a_2z^2+\cdots\eqno{(1.1)}$$
which are analytic and univalent in the unit disk $E=\{z\colon
|z|<1\}$ and map $E$ onto some domains $D$. A function $f\in S$ is said to be starlike if the domain $D$ is starlike with respect to the origin. The family of starlike functions is denoted by $S^\ast$. It is known that a function $f\in S$ is starlike if and only if it satisfies the geometric condition:
$$Re\frac{zf'(z)}{f(z)}>0.\eqno{(1.2)}$$

The above geometric condition implies that the quantity $zf'(z)/f(z)$ belongs to the class $P$ of analytic functions:
$$p(z)=1+c_1z+c_2z^2+\cdots\eqno{(1.3)}$$
which have positive real part in $E$. The family of starlike functions in the unit disk has attracted much attention in the past. The volume of work being published on this family of functions leaves no one in doubt about the importance attached to it both in the past and the present.

In this paper we determine, under certain conditions, the best possible upper bounds on some functionals defined in the coefficients space of starlike functions. These are functionals which have arisen from the study of coefficient problems of certain family of univalent functions. The study of functionals similar to those being considered in this paper is not new. For instance, bounds on $|a_n|$ and $|a_3-\lambda a_2^2|$ can be found in many literatures (see for example \cite{PN}-\cite{RS}). In particular, the functional $|a_3-\lambda a_2^2|$ is known as the Fekete-Szego functional for both real and complex values of the parameter $\lambda$. The problem of determining the sharp bound on the Fekete-Szego functional has its origin in a conjecture of Littlewood and Parley (1932) that the true bound on the coefficients of an odd univalent function is 1, which was disproved in 1933 be Fekete and Szego via the determination of the sharp bound on the functional (see \cite{PL}).

The functional has since continued to recieve attention of researchers in geometric function theory.

In Section 3, we consider functionals of the form $|a_4-\gamma a_2a_3|$, $|a_4-\gamma a_2a_3-\eta a_2^3|$, $|a_5-\mu a_2^2a_3|$ and $|a_5-\xi a_2a_4-\zeta a_2^3|$ where the parameters $\gamma, \eta, \mu, \xi, \zeta$ are all real numbers. These functionals have been found to have applications in certain coefficient problems, which are of independent interest.
 \medskip

{\bf 2.0 Preliminary Lemmas}\vskip 2mm

We shall need the following well known inequalities.\vskip 2mm

{\bf Lemma 2.1}(\cite{PN}-\cite{RS})\vskip 2mm

{\em Let $p\in P$. Then $|c_k|\leq 2$, $k=1,2,3,\cdots$. Equality is attained for the Moebius function
$$L_0(z)=\frac{1+z}{1-z}.\eqno{(2.1)}$$}

{\bf Lemma 2.2}(\cite{PL,SK})\vskip 2mm

{\em Let $p\in P$. Then $$\left|c_2-\frac{c_1^2}{2}\right|\leq 2-\frac{|c_1|^2}{2}\eqno{(2.2)}$$
The result is sharp. Equality holds for the function
$$p(z)=\frac{1+\frac{1}{2}(c_1+\varepsilon\bar{c_1})z+\varepsilon z^2}{1-\frac{1}{2}(c_1-\varepsilon\bar{c_1})z-\varepsilon z^2},\;\;|\varepsilon|=1.\eqno{(2.3)}$$}

Note that the inequality (2.2) can be written as
$$c_2=\frac{1}{2}c_1^2+\varepsilon\left(2-\frac{1}{2}|c_1|^2\right),\;\;|\varepsilon|\leq 1.\eqno{(2.4)}$$

\medskip

{\bf 3.0 Main Result}\vskip 2mm

{\bf Theorem 3.1}\vskip 2mm

{\em Let $f(z)$ given by (1.1) be starlike function. Then for real numbers $\gamma, \eta, \mu, \xi, \zeta$ such that $1-\gamma, 1-2\mu, 1-\xi, 1-2\zeta$ and $1-2\xi-2\zeta$ are all nonnegative, we have the sharp inequalities:

$$|a_4-\gamma a_2a_3|\leq 4-6\gamma;\;\;if\;\;\gamma\leq\frac{5}{9},$$

$$|a_4-\gamma a_2a_3-\eta a_2^3|\leq
4-6\gamma-8\eta;\;\;if\;\;3\gamma+4\eta\leq\frac{5}{3},$$

$$|a_5-\mu a_2^2a_3|\leq 5-12\mu;\;\;if\;\;\mu\leq
\frac{2}{9},$$

$$|a_5-\xi a_2a_4-\zeta a_3^2|\leq
5-8\xi-9\zeta;\;\;if\;\;5\tau+9\omega\leq 2.$$}

\begin{proof} Since $f(z)$ is starlike, there exists $p\in P$ such that
$$zf'(z)=p(z)f(z)\eqno{(3.1)}$$
Comparing coefficients of both sides of (3.1) using (1.1) and (1.3) we see that
$$a_2=c_1$$
$$2a_3=c_2+c_1^2$$
$$6a_4=2c_3+3c_2c_1+c_1^2$$
$$24a_5=6c_4+8c_3c_1+6c_2c_1^2+3c_2^2+c_1^4$$
$$\cdots\;\cdots\;\cdots\;\cdots\;\cdots$$
so that
$$a_4-\gamma a_2a_3=\frac{c_3}{3}+(1-\gamma)\frac{c_1}{2}\left\{c_2+\frac{2(1-3\gamma)}{3(1-\gamma)}\frac{c_1^2}{2}\right\}\eqno{(3.2)}$$
$$a_4-\gamma a_2a_3-\eta a_2^3=\frac{c_3}{3}+(1-\gamma)\frac{c_1}{2}\left\{c_2+\frac{2(1-3\gamma-6\eta)}{3(1-\gamma)}\frac{c_1^2}{2}\right\}\eqno{(3.3)}$$
$$a_5-\mu a_2^2a_3=\frac{c_4}{4}+\frac{c_3c_1}{3}+\frac{c_2^2}{8}+(1-2\mu)\frac{c_1^2}{4}\left\{c_2+\frac{1-12\mu}{3(1-2\mu)}\frac{c_1^2}{2}\right\}\eqno{(3.4)}$$
$$\aligned
a_5-\xi a_2a_4-\zeta a_3^2
&=\frac{c_4}{4}+(1-\xi)\frac{c_3c_1}{3}+(1-2\zeta)\frac{c_2^2}{8}\\
&+(1-2\xi-2\zeta)\frac{c_1^2}{4}\left\{c_2+\frac{1-4\xi-6\zeta}{3(1-2\xi-2\zeta)}\frac{c_1^2}{2}\right\}\endaligned\eqno{(3.5)}$$

Recall that the real numbers $1-\gamma, 1-2\mu, 1-\xi, 1-2\zeta$ and $1-2\xi-2\zeta$ are nonnegative. We eliminate $c_2$ in each of the terms in the curly brackets in (3.2) - (3.5) using the equality (2.4). For instance, we have from (3.2),
$$c_2+\frac{2(1-3\gamma)}{(3(1-\gamma)}\frac{c_1^2}{2}=\frac{5-9\gamma}{3(1-\gamma)}\frac{c_1^2}{2}+\varepsilon\left(2-\frac{|c_1|^2}{2}\right).\eqno{(3.6)}$$

Since $2-\frac{|c_1|^2}{2}\geq 0$, the absolute value of (3.6) attains its maximum for $|c_1|=2$ provided $\gamma\leq\frac{5}{9}$ (which is the condition given in the first inequality of the theorem). Thus (3.6) yields  
$$\left|c_2+\frac{2(1-3\gamma)}{3(1-\gamma)}\frac{c_1^2}{2}\right|\leq\frac{2(5-9\gamma)}{3(1-\gamma)},\eqno{(3.7)}$$
so that, by triangle inequality and Lemma 2.1, (3.2) yields the first inequality of the theorem. Similar arguments and computations from (3.3) to (3.5), lead to the remaining inequalities respectively.

For each of the real numbers $\gamma, \eta, \mu, \xi$ and $\zeta$, equality is attained in each case by the Koebe function (up to rotations) given by:
$$k(z)=\frac{z}{(1-z)^2}.\eqno{(3.8)}$$
\end{proof}
 
\medskip

{\bf 4.0 Conclusion}\vskip 2mm

The study of functionals in the theory of analytic and univalent functions is here boosted with the consideration of new ones. The functionals considered in this work are closely associated with certain coefficient problems in geometric functions theory.

\bigskip

\bibliographystyle{amsplain}

\end{document}